\documentclass[12pt]{amsart}
\title[Quasiregular maps of Sierpi\'nski carpet Julia sets]{Quasiregular maps of 
Sierpi\'nski carpet\\ Julia sets}
\author{Sergei Merenkov, Letian Shen}

\address{Sergei Merenkov\\Department of Mathematics, The City College of New York and CUNY Graduate Center, New York, NY 10031, USA.}
\email{smerenkov@ccny.cuny.edu}
\thanks{SM supported by NSF grant DMS-2247364.}

\address{Letian Shen\\Department of Mathematics\\University of Science and Technology of China, Anhui 230026\\China}
\email{ltshen@mail.ustc.edu.cn}

\usepackage[matrix,arrow,curve,cmtip]{xy}
\usepackage{amsmath, amssymb, amsthm,latexsym}
\usepackage{graphicx}

\newtheorem{thm}{Theorem}[section]

\newtheorem{lem}[thm]{Lemma}
\newtheorem{prop}[thm]{Proposition}

\newtheorem{question}{Question}[section]

\newtheorem*{ex}{Example}
\newtheorem*{pf}{Proof}

\newcommand{\R}{\mathbb{R}}
\newcommand{\C}{\mathbb{C}}
\newcommand{\D}{\mathbb{D}}  

\newcommand{\N}{\mathbb{N}}

\newcommand{\ra}{\rightarrow}
\newcommand{\hC}{\widehat{\C}}
\newcommand{\sub}{\subseteq }
\newcommand{\J}{\mathcal{J}}

\newcommand{\co}{\colon}

\makeatletter
\renewenvironment{pf}[1][\proofname]{%
  \par
  \pushQED{\qed}%
  \normalfont\topsep6\p@ \@plus6\p@ \trivlist
  \item[\hskip\labelsep\itshape
        \bfseries #1\@addpunct{:}]%
  \ignorespaces
}{%
  \popQED\endtrivlist
}
\makeatother

\makeatletter
\renewenvironment{ex}[1][\examplename]{%
  \par
  \pushQED{\qed}%
  \normalfont\topsep6\p@ \@plus6\p@ \trivlist
  \item[\hskip\labelsep\itshape
        \bfseries #1\@addpunct{:}]%
  \ignorespaces
}{%
  \popQED\endtrivlist
}

\renewenvironment{ex}{%
  \normalfont\topsep6\p@ \@plus6\p@ \trivlist
  \item[\hskip\labelsep\bfseries Example:]%
  \ignorespaces
}{%
  \endtrivlist
}

\makeatother

\DeclareMathOperator*{\esssup}{ess\,sup}

\begin{document}

\abstract{We prove that if $f$ and $g$ are postcritically finite rational maps whose Julia sets $\J(f), \J(g)$, respectively, are Sierpi\'nski carpets, and if $\xi$ is a quasiregular map of the Riemann sphere $\hC$ with $\xi^{-1}(\J(g))=\J(f)$, then $\xi$ is the restriction of a rational map to the Julia set $\J(f)$. Moreover, when $g=f$ we prove that, for some positive integers $k$ and $l$, $f^k\circ \xi^l=f^{2k}$. These conclusions extend the main results of \cite{BLM}.} Finally, we demonstrate that when Julia sets of postcritically finite rational maps are not Sierpi\'nski carpets, say they are tree-like or gaskets, the above conclusions no longer hold.
\endabstract

\maketitle

\section{Introduction}
\noindent A \textit{Sierpi\'nski carpet}, or \emph{carpet} for short, is a homeomorphic image of the standard Sierpi\'nski carpet fractal. 
Due to~\cite{Wh}, a subset $S$ of $\hC$ is a Sierpi\'nski carpet if and only if it has empty interior and can be written as 
$$S=\hC\backslash \bigcup_{k\in \N}D_k,$$
where $D_k,\, k\in\N$, are pairwise disjoint Jordan domains with 
$\partial D_k\cap \partial D_l=\emptyset$ for $k\neq l$, and $\mathrm{diam}(D_k)\to 0$ as $k\to \infty$. The boundaries $\partial D_k$ of $D_k,\ k\in\N$, are called \emph{peripheral circles}. A carpet $S$ is called a \textit{round Sierpi\'nski carpet}, or \textit{round carpet}, if the domains $D_k,\, k\in \N$, are open geometric disks.

A rational map $f$ is called \emph{postcritically finite} if each cri\-ti\-cal point of $f$ has a finite forward orbit.
The following result is Theorem~1.4 from~\cite{BLM}.
\begin{thm} [\cite{BLM}]
	Let $f$ and $g$ be two {postcritically finite} rational maps and suppose that the corresponding Julia sets $\J(f)$ and $\J(g)$ are Sierpi\'nski carpets.
If $\xi$ is a quasisymmetric homeomorphism of $\J(f)$ onto $\J(g)$,
 then it is the restriction to $\J(f)$ of a Mobius transformation on $\hC$.
\end{thm}


Our first main result extends this theorem to quasiregular maps.
\begin{thm}\label{thm:main1}
	Let $f$ and $g$ be postcritically finite rational maps whose Julia sets $\J(f), \J(g)$, respectively, are Sierpi\'nski carpets.
 If $\xi\co\hC\to \hC$ is a quasiregular map such that
 $\xi^{-1}(\J(g))=\J(f)$, 
 then $\xi|_{\J(f)}$ has a unique extension to a rational map $\widehat{\xi}\co\hC\ra \hC$ satisfying
 \begin{equation}\label{eqn:fxi1}
 	G^{n+k}\circ\widehat\xi=G^n\circ\widehat\xi\circ F^k,
 \end{equation} 
 for all $k,n\in\N, \ n\ge l$, for some $l\in\N$, where $F$ and $G$ are some iterates of $f$ and $g$, respectively. 
\end{thm}
\noindent Moreover, when $f=g$, we have the following result.
\begin{thm}\label{thm:main2}
  Let $f$ be a postcritically finite rational map whose Julia set is a Sierpi\'nski carpet, let $\xi\co\hC\to \hC$ be a quasiregular map such that
  $\xi^{-1}(\J(g))=\J(f)$, and let $\widehat\xi$ be the unique rational extension of $\xi|_{\J(f)}$ from Theorem~\ref{thm:main1}. If $\mathrm{deg}(\widehat{\xi})>1$, then $\widehat{\xi}$ is postcritically finite, the Julia set $\J(\widehat\xi)$ of $\widehat\xi$ equals $\J(f)$,
  and there exist $k, l\in \N$ such that
  \begin{equation}\label{eqn:fxi2}
  f^k\circ (\widehat{\xi})^l=f^{2k}.
  \end{equation}
\end{thm}

The structure of this paper is as follows. Section \ref{sec:qr} introduces  qua\-si\-re\-gu\-lar and related maps. Section \ref{sec:Sch} is dedicated to the theory of Schottky maps, as developed by the first author~\cite{Me1,Me2, Me3}, where the convergence Theorem~\ref{thm:Der} is an important ingredient of Proposition~\ref{prop:basic} in Section~\ref{sec:review}. 
Section \ref{sec:conformalelevator} discusses the technique of conformal elevator from~\cite[Section 4]{BLM}, which enables us to construct local quasiconformal maps with images of definite size. This is another ingredient needed for Proposition~\ref{prop:basic}. 
Section \ref{sec:review} outlines the proof of Proposition~\ref{prop:basic} adapted from \cite{BLM} to a quasiregular case.
After some preparation in Section \ref{sec:equilibrium} on rational maps sharing a measure of maximal entropy, based on~\cite{Le, Le2, LP, Ye}, we prove the main Theorems \ref{thm:main1} and \ref{thm:main2} in Section \ref{sec:main}. 
Section \ref{sec:other} discusses quasiregular maps between Julia sets of topologies other than the Sierpi\'nski carpet. This section is inspired by the work of Y.~Luo and D.~Ntalampekos~\cite{LN}, which reveals that gaskets do not exhibit the local rigidity properties of carpets. Our examples in this section add evidence to the contrast between carpet and gasket as well as tree-like Julia sets of polynomials.  

\medskip
\noindent
{\bf Acknowledgements.} The authors thank the Institute for Ma\-the\-ma\-ti\-cal Sciences at Stony Brook University for its hospitality.
The second author thanks the University of Science and Technology of China for the support of the travel that resulted in this work. The authors also thank Yusheng Luo for useful comments and references.

\section{Quasiregular maps and related concepts}\label{sec:qr}
\noindent A non-constant mapping $f\co U\to \hC$ on a domain $U\sub \hC$ is called \textit{$K$-quasiregular}, $K\geq 1$, if $f$ is in the Sobolev space $W_\mathrm{loc}^{1,2}$ and for almost every $z\in U$,
$$
  \left\lVert Df(z)\right\rVert^2\leq K\mathrm{det}(Df(z)).
$$
If $f$ is $K$-quasiregular for some $K\geq 1$, then $f$ is called \textit{quasiregular}. An equivalent condition for $f\co U\to \C$ to be quasiregular is that $f\in W_\mathrm{loc}^{1,2}$ and $f$ satisfies the \textit{Beltrami equation}
$$
\frac{\partial f}{\partial \bar{z}}=\mu(z)\frac{\partial f}{\partial z}
$$
for some complex valued Lebesgue measurable function $\mu$ with 
$$
\esssup|\mu|<1.
$$
If a $K$-quasiregular map $f$ is a homeomorphism onto its image in addition, then $f$ is called a $K$-\textit{quasiconformal} map. A 1-quasiconformal map is called \emph{conformal}.

We note here a useful factorization of quasiregular maps by S.~Sto\-\"i\-low~\cite{St} (see also~\cite[Section~5.5]{AIM} and~\cite{LP2}):
\begin{thm}[Sto\"ilow's Factorization]
  If $f\co U\to \hC$ is a quasiregular map, then there exists a homeomorphism $h\co U\to V$ which satisfies the same Beltrami
  equation as $f$. Moreover, for every such $h$ there exists a holomorphic map $g\co V\to \hC$ such that 
  $$f=g\circ h.$$
\end{thm} 
\noindent It follows easily from this theorem that a quasiregular map $f\co U\to \hC$ is a local homeomorphism at all but a discrete set of points.
We call a point $p\in U$ \textit{critical} if $f$ is not a local homeomorphism in any neighborhood of $p$.  The set of critical points of $f$ will be denoted by ${\rm crit }(f)$.

We will also need the following quasiregular version of Montel’s theorem~\cite[Corollary 5.5.7]{AIM} as well as a compactness property of quasisymmetric maps below. 
\begin{thm}[\cite{AIM}]\label{thm:converge1}
Suppose $z_0, z_1\in\C$ and let $h_k\colon\Omega\to \C, k = 1, 2,\dots$, be a
sequence of $K$-quasiregular mappings defined on a domain $\Omega\subset \C$, each omitting
these two values. Then there is a subsequence converging locally uniformly on $\Omega$
to a mapping $h$,
and $h$ is either a $K$-quasiregular mapping or a constant.
\end{thm}

Finally, recall that if $(X, d_X)$ and $(Y, d_Y)$ are two metric spaces, a homeomorphism $f\co X \to Y$ is said to be $\eta$-\emph{quasisymmetric} if there is an increasing homeomorphism $\eta\co [0, \infty) \to [0, \infty)$ such that for any triple $x, y, z$ of distinct points in $X$, we have
$$\displaystyle {\frac {d_{Y}(f(x),f(y))}{d_{Y}(f(x),f(z))}}\leq \eta \left({\frac {d_{X}(x,y)}{d_{X}(x,z)}}\right).$$

\section{Schottky maps and their local rigidity}\label{sec:Sch}
\noindent A \textit{relative Schottky set} $S$ in a domain $D\sub \hC$ is a set of the form 
$$S=D\backslash \cup_{i\in I}B_i,$$
where each $B_i,\, i\in I$, is an open geometric disk with the closure $\overline{B_i}$ contained in $D$, and $\overline{B_i}\cap\overline{B_j}=\emptyset$ for $i\neq j$.
The boundary $\partial B_i$ of a disk $B_i,\, i\in I$, is called a \textit{peripheral circle}. If $D=\C$  or $\hC$, $S$ is called a \textit{Schottky set}. The study of Schottky sets was initiated in~\cite{BKM}, and relative Schottky sets and maps in~\cite{Me1, Me2}.


Let $S\sub \C$ be an arbitrary subset that has an accumulation point $p\in S$, and let $U$ be a neighborhood of $p$ in $\C$.
Suppose that $f\co S\cap U\to \C$ is a continuous function. We say that $f$ is \textit{conformal} at $p$ if 
\begin{equation}\label{E:Conformal}
	\lim_{q\in S,\, q\to p}\frac{f(q)-f(p)}{q-p}
\end{equation}
exists and is non-zero. The limit above is called the \textit{derivative} of $f$ at $p$ and is denoted by $f'(p)$.

Suppose that $S=D\backslash \cup_{i\in I}B_i$ and $\widetilde{S}=\widetilde{D}\backslash \cup_{j\in J}\widetilde{B_i}$ are relative Schottky sets.
Let $U$ be an open subset of $D$ and $f\co S\cap U\to \widetilde{S}$ is a local homeomorphism.
It follows that $p\in S\cap U\cap \partial B_i$ for some $i\in I$ if and only if 
$f(p)\in \partial \widetilde{B_j}$ for some $j\in J$. Such a map $f\co S\cap U\to \widetilde{S}$ is called a \textit{Schottky map} if it is conformal at every point of $S\cap U$, and the derivative $f'$ is a continuous function on $S\cap U$. We use standard modifications when either $p$ of $f(p)$ is $\infty$ in Equation~\eqref{E:Conformal} above.

The notion of Schottky maps is closely related to that of quasiconformal maps. In fact, we will use the following criterion \cite[Lemma~2.2]{Me3} for maps to be Schottky maps. 
\begin{lem}[\cite{Me3}]\label{lem:Schmaps} 
  Let $S\sub \hC$ be a Schottky set of measure zero. 
Suppose   $U\sub \hC$ is  open and   $f\co  U\ra \hC$ is  a quasiconformal map with $f^{-1}(S)=S\cap U$.
Then $f\co S\cap U\to S$ is a Schottky map. 
In particular, if $f\co \hC\ra \hC$ is a quasiregular map with $f^{-1}(S)=S$, then 
$f\co S\setminus \mathrm{crit}(f)\to S$ is a Schottky map. 
\end{lem}

An important property of Schottky maps is their $local\;rigidity$, i.e., they are uniquely determined by their local behavior near one point. We need the following notion to state related results, Theorem~\ref{thm:Un}~\cite[Corollary~4.2]{Me2} and Theorem~\ref{thm:Der}~\cite[Theorem~5.2]{Me2} below. A relative Schottky set $S=D\setminus\cup_{i\in I} B_i$ is called \emph{locally porous at} $p\in S$ if there exist a neighborhood $U$ of $p$ and constants $r_0>0,\, C\geq 1$, such that for every $q\in S\cap U$ and each $r$ with $0<r\leq r_0$, there exists $i\in I$ with $B(q,r)\cap B_i\neq \emptyset$ and
$$
r/C\leq {\rm diam}(B_i)\leq Cr,
$$
where $B(q,r)$ denotes the open disk of radius $r$ centered at $q$, and ${\rm diam}(X)$ denotes the diameter of a set $X\subseteq\C$. A relative Schottky set $S$ is called \emph{locally porous} if it is locally porous at every $p\in S$. Every locally porous relative Schottky set has measure zero since it cannot have Lebesgue density points.

\begin{thm}[\cite{Me2}]\label{thm:Un}
Let $S$ be a locally porous relative Schottky set in $D\subseteq \C$,  and suppose that $U\subseteq D$ is an open  set such that $S\cap U$ is connected. 
Let $f$ and $g\co  S\cap U\to \widetilde S$ be Schottky maps into a relative Schottky set $\widetilde S$ in a domain $\widetilde D$, and consider 
$$
E=\{p\in S\cap U\co f(p)=g(p)\}.
$$  
Then $E=S\cap U$, provided $E$ has an accumulation point in $U$. 
\end{thm} 
We also record another theorem \cite[Theorem 5.2]{Me2} on the rigidity of Schottky maps. Roughly speaking, it says that, under a mild assumption that is satisfied for any subhyperbolic dynamical system, a sequence $f_k$ of locally uniformly converging Schottky maps will become constant when $k$ is large enough. 
\begin{thm}[\cite{Me2}]\label{thm:Der}
Let $S$ be a locally porous relative Schottky set in a  domain $D\subseteq\C$, and $p\in S$ be an arbitrary point. 
Suppose that $U\subseteq D$ is an open neighborhood of $p$ such that $S\cap U$ is connected. We assume that there exists a Schottky map 
$g\co S\cap U\to S$ with $g(p)=p$ and $g'(p)\neq1$.
Let $\widetilde S$ be a relative Schottky set in a domain $\widetilde D$ and let $(f_k)_{k\in\N}$ be a sequence of Schottky maps $f_k\co S\cap U\to \widetilde S$.  We assume that for each $k\in\N$ there exists an open set $\widetilde U_k$ so that the map $f_k\co S\cap U\to \widetilde S\cap \widetilde U_k$ is a homeomorphism, and the sequence $(f_k)$
converges locally uniformly to a  homeomorphism $f\co S\cap U\to \widetilde S\cap\widetilde U$, where $\widetilde U$ is an open set. 
Then there exists $N\in \N$ such that $f_k=f$ in $S\cap U$ for all $k\geq N$.
\end{thm}

\section{Conformal elevator}\label{sec:conformalelevator}
\noindent In \cite[Section 4]{BLM}, the authors introduced the technique of \textit{conformal elevator} for subhyperbolic rational maps. Here we discuss this technique for postcritically finite maps
only, and we will apply it to the proof of Proposition~\ref{prop:basic} in Section~\ref{sec:review}.

Assume that $f$ is a postcritically finite rational map with 
$$\J(f)\sub \frac{1}{2}\D\quad and \quad f^{-1}(\D)\sub \D.$$
For a map $f$, say with an attracting fixed point, these are easy to achieve by conjugating $f$ with an appropriate M\"obius transformation. 
In this subsection we use Euclidean metric on $\C$.
Fix $\epsilon_0>0$ so small that $\mathrm{diam}(\J(f))>2\epsilon_0$, and so that every disk $B'=B(q,r')$ centered at a point $q\in \J(f)$ with radius
$r'\leq 8\epsilon_0$ is contained in $\D$, and contains at most one point in the postcritical set of $f$, i.e., in $\mathrm{post}(f)=\cup_{n\in\N}f^n({\rm crit}(f))$.   
Let $B=B(p,r)$ be a small disk centered at a point $p\in \J(f)$ and of radius $r<\epsilon_0$. We know that $\J(f)\sub f^n(B)$ for sufficiently large 
$n$ (see \cite[Theorem 4.2.5 (ii)]{Be}), thus $\mathrm{diam}(f^n(B))>2\epsilon_0$. 
In other words, there exists a maximal number $n\in \N$ such that $\mathrm{diam}(f^n(B))\leq 2\epsilon_0$.
By uniform continuity of $f$ near the compact set $\J(f)$, there exists $\delta_0>0$ only depending on $f$ such that 
$$\mathrm{diam}(f^n(B))\geq \delta_0.$$
Moreover, we have a lower bound for diameters of image sets $f^n(A)$, where $A\sub B$ is connected (\cite[Lemma 4.1]{BLM}):
\begin{lem}[\cite{BLM}]\label{lem:conformalelevator}
There exist constants $\gamma>0$ and $C>0$ only depending on $f$ and independent of $B=B(p,r)$, such that for any connected set $A\sub B$,
$$\mathrm{diam}(f^n(A))^\gamma\geq C\frac{\mathrm{diam}(A)}{\mathrm{diam}(B)}.$$
\end{lem}
\begin{pf}
We sketch the proof here. First, by our choice of $\epsilon_0$, 
we find a disk $B'=B(q,r')$ with $q\in \J(f)$ and $r'\leq 8\epsilon_0$, such that 
$$f^n(B)\sub \frac{1}{2}B'\sub \D \quad {\rm and} \quad B'\cap \mathrm{post}(f)\sub \{q\}.$$
Let $\Omega\in \hC$ be the unique component of $f^{-n}(B')$ that contains $B$.\\
Second, there exist conformal maps $\varphi\co B'\to \D$ and $\psi\co \Omega\to \D$ such that 
$$(\varphi\circ f^n\circ \psi^{-1})(z)=z^k$$
for all $z\in \D$, where $k\leq N=N(f)$ is uniformly bounded only by $f$ and independent of $B$. 
Third, 
by Koebe distortion theorem we have 
$$|u-v|\approx  \left\lvert (\psi^{-1})'(0)\right\rvert \cdot \left\lvert \varphi(u)-\varphi(v)\right\rvert, $$
whenever $u,v\in B$, and where we use $a\approx b$ to indicate that there is a constant $C\ge1$ independent of the disk $B$ such that $a/C\le b\le Cb$. In particular, 
$$\mathrm{diam}(B)\approx  \left\lvert (\psi^{-1})'(0)\right\rvert \cdot \mathrm{diam}(f^n(B)).$$
On the other hand, 
\begin{align*}
  1\approx  \mathrm{diam}(f^n(B))  &\approx  \mathrm{diam}(\varphi(f^n(B))) 
  \\& =\mathrm{diam}(P_k(\psi(B)))\le C\cdot \mathrm{diam}(\psi(B))\leq 2,
\end{align*}
where $C>0$ is a constant independent of $B$.
Hence $\mathrm{diam}(\psi(B))\approx  1$, and so $\mathrm{diam}(B)\approx \left\lvert (\psi^{-1})'(0)\right\rvert $.
This implies that 
$$\frac{|u-v|}{\mathrm{diam}(B)}\approx  \left\lvert \varphi(u)-\varphi(v)\right\rvert,$$
whenever $u,v\in B$. Hence
\begin{align*}
\frac{\mathrm{diam}(A)}{\mathrm{diam}(B)}\approx  \mathrm{diam}(\psi(A)) &\le C\cdot  \mathrm{diam}(P_k(\psi(A)))^{1/N}\\&=\mathrm{diam}(\varphi(f^n(A)))^{1/N}\approx  \mathrm{diam}(f^n(A))^{1/N},
\end{align*}
for another constant $C>0$ independent of $B$. 
\end{pf}

\section{Local quasiconformal maps between carpet Julia sets}\label{sec:review}
\noindent In this section, $f$ and $g$ are postcritically finite rational maps, and the corresponding Julia sets $\J(f)$ and $\J(g)$ are Sierpi\'nski carpets.
 In \cite[Theorem~1.4]{BLM}, the authors proved that a quasisymmetric map $\xi\co \J(f)\to\J(g)$ is the restriction to $\J(f)$ of a M\"obius transformation. 
 Many of the arguments in that paper are local and it makes sense to investigate whether they carry over to {quasiregular} maps $\xi$ that restrict to maps between $\J(f)$ and $\J(g)$. Note that a quasisymmetric map $\xi\co\J(f)\to\J(g)$ extends to a global quasiconformal map of $\hC$. One of the main ingredients of~\cite{BLM} is the functional equation~\cite[Equation~(1.2)]{BLM} 
 $$
  g^{m'}\circ \xi = g^m\circ \xi \circ f^n,
 $$
established on all of $\J(f)$ in~\cite[Section~8, Step IV]{BLM}.
Since quasiregular maps are locally quasiconformal away from critical points, 
we obtain such a functional equation in the local setting; see Proposition~\ref{prop:basic} below, and then conclude that it holds globally on $\J(f)$ for quasiregular maps in Theorem~\ref{thm:basicthm}. 
We outline the proofs with these modifications for readers' convenience. 
\begin{prop}\label{prop:basic}
Let $f, g\co \hC \to \hC$ be postcritically finite rational maps, and suppose that the corresponding Julia sets $\J(f)$ and $\J(g)$ are Sierpi\'nski carpets. 
If $U$ and $V$ are open sets non-trivially intersecting $\J(f)$ and $\J(g)$, respectively, and $\xi\co U \ra V$ is a quasiconformal map with $\xi^{-1}(\J(g)\cap V)=\J(f)\cap U$, then there exist an open set $W\subseteq U$ with $\J(f)\cap W\neq\emptyset$, and $n, m, m'\in \N$,
such that on $\J(f)\cap W$ one has
\begin{equation}\label{eq:keyeq}
  g^{m'}\circ \xi = g^m\circ \xi \circ f^n.
\end{equation}
\end{prop}

\begin{pf}
  By conjugating $g$ (hence $\J(g)$), we may assume that $\J(g)\subset \D$. Below we use the Euclidean metric.
  Since repelling periodic points of $f$ are dense in $\J(f)$ (see
 \cite[p.~148, Theorem~6.9.2]{Be}), we can find a repelling periodic point 
$p\in \J(f)\cap U$ of $f$ of period $n$ and $\rho>0$ such that $U_0=B(p, \rho)\sub U$  is disjoint from $\mathrm{post}(f)$.  
 Denote $n\cdot k$ by $n_k$, and the connected component of $f^{-n_k}(U_0)$ containing $p$ by $U_k$. 
 By possibly taking a smaller $\rho$, we may assume that $U_0\supseteq  U_1\supseteq  U_2\supseteq ...$ and $\mathrm{diam}(U_k)\to 0$ as $k\to\infty$. 
 We take the branch $f^{-n_k}$ on $U_0$ to be the inverse map of 
 $f^{n_k}|_{U_k}\co U_k\to U_0$. Then the choice of these inverse branches is {\em consistent} in the sense that we have 
 \begin{equation}\label{eqn:consistent}
 f^{n_{k+1}-n_k}\circ f^{-n_{k+1}}=f^{-n_k}
 \end{equation}
 for all $k$. 
 Let $r_k>0$ be the smallest number such that $\xi(U_k)\sub B(\xi(p),r_k)$. Since $\xi(p)\in \xi(U_k)$, $\mathrm{diam}(\xi(U_k))\approx r_k$, where, similar to above,  the symbol $\approx$ between two quantities means that their quotient is bounded above and below by a constant independent of $k$. 
By applying Lemma~\ref{lem:conformalelevator} to $A=\xi(U_k)$ and $B=B(\xi(p),r_k)$, we can find iterates
$g^{m_k}$ such that $\mathrm{diam}\left(g^{m_k}(\xi(U_k))\right)\approx 1$. Since $g^{m_k}(\xi(U_k))$ meets $\J(g)$ at $g^{m_k}(\xi(p))$, it is contained in a 
bounded region of $\C$.

Now we consider the composition 
$$
h_k=g^{m_k}\circ \xi\circ f^{-n_k}
$$ 
defined on $B(p,\rho)$ for $k\in \N$. 
Since $h_k$ are uniformly $K$-quasiregular ($K$ only depends on $\xi$) maps from $U_0$ into a compact subset of $\C$,  
Theorem~\ref{thm:converge1} implies that a subsequence of $\{h_k\}$
subconverges locally uniformly to a map $h\co U_0\to \C$  that 
is also quasiregular or a constant. The latter possibility is ruled out since 
$\mathrm{diam}\left(h_k(U_0)\right)=g^{m_k}(\xi(U_k))\approx 1$. For simplicity of notations we assume that the whole sequence $h_k$ converges; though $n_k$ is no longer $n\cdot k$ in this case. The map
$h$ has at most countably many critical points, and so there exists a non-critical point $q\in \J(f)\cap U_0$ and a small neighborhood $B(q,r)$ of $q$ such that $h$ is injective on $B(q,r)$, and hence quasiconformal. 
In what follows, we replace $U_0$ by this $B(q,r)$.

We know by \cite[Theorem 1.10]{BLM} that the peripheral circles of a Sierpi\'nski carpet Julia set $\J(f)$ of a postcritically finite rational map $f$ satisfy nice geometric properties. Namely, they are uniformly relatively separated uniform quasicircles and appear at all locations and scales; see~\cite{BLM} for definitions. In particular, 
according to \cite[Corollary 1.2]{Bo}, there exists a quasisymmetric map $\beta$ on $\hC$ such that 
$S=\beta(J(f))$ is a round Sierpi\'nski carpet, i.e., its peripheral circles are geometric circles. Moreover, $S$ is locally porous.
%

From Lemma~\ref{lem:Schmaps} we know that any iterate 
$$f_\beta^n=\beta\circ f^n\circ \beta^{-1}\co S\backslash \mathrm{crit}(f_\beta^n)\to S$$ is a 
Schottky map.  
Moreover, since $$h^{-1}(\J(g))=h_k^{-1}(\J(g))=\J(f)\cap U_0,$$
by the same lemma, the quasiconformal maps 
$\beta\circ h\circ \beta^{-1}, \beta\circ h_k\circ \beta^{-1},\ k\in\N$, 
are Schottky maps on $S\cap \beta(U_0)$.

It is easy to see (\cite[Lemma~2.1]{Me3}) that one can find a small Jordan region $W'$ with $a\in W'\sub \beta(U_0)$
such that $S\cap W'$ is connected. 
Since repelling periodic points are dense in $\J(f)$, we can find 
a point $a\in \left(S\backslash \mathrm{crit}(f_\beta^n)\right)\cap W'$ and $n\in\N$  such that 
$$f_\beta^n(a)=a \quad \mathrm{and}\quad (f_\beta^n)'(a)\neq 1.$$
Theorem~\ref{thm:Der} implies that 
$$\beta\circ h_k\circ \beta^{-1}\equiv \beta\circ h\circ \beta^{-1}\quad  \mathrm{on}\quad S\cap W', $$
 for $k$ large enough. Hence
$h_k=h_l$ on $\J(f)\cap W''$, where $W''=\beta^{-1}(W')$, for some $k<l$, i.e., 
$$g^{m_k}\circ \xi\circ f^{-n_k}=g^{m_l}\circ \xi\circ f^{-n_l}\quad  \mathrm{on}\quad \J(f)\cap W''.$$ 
Finally, we pre-compose $f^{n_l}$ on both sides and from Equation~\eqref{eqn:consistent} get 
$$g^{m_k}\circ \xi\circ f^{n_l-n_k}=g^{m_l}\circ \xi \quad  \mathrm{on}\quad \J(f)\cap W,$$ where $W=f^{-n_l}(W'')$, establishing~\eqref{eq:keyeq}.
\end{pf}

As a consequence we obtain the following result.
\begin{thm}\label{thm:basicthm}
	If $f$ and $g$ are postcritically finite rational maps whose respective Julia sets $\J(f), \J(g)$ are Sierpi\'nski carpets, and 
  if $\xi$ is a quasiregular map $\xi\co \hC\to \hC$ such that $\xi^{-1}(\J(g))=\J(f)$,
  then Equation~\eqref{eq:keyeq} holds everywhere on $\J(f)$.
\end{thm}

\begin{pf}
	Since $\J(g)$ has the same properties as $\J(f)$, by \cite[Corollary 1.2]{Bo} there exists a quasisymmetric map $\gamma$ on $\hC$ such that 
	$\widetilde S=\beta(J(g))$ is a round Sierpi\'nski carpet.
  Using the same notations and arguments as in Proposition~\ref{prop:basic} and its proof, we know that $H_1=\gamma\circ g^{m'}\circ \xi\circ \beta^{-1}$ and $H_2=\gamma\circ g^m\circ \xi \circ f^n\circ \beta^{-1}$
are Schottky maps from $S\backslash (\mathrm{crit}(H_1)\cup \mathrm{crit}(H_2))$ to $\widetilde S$. 
Proposition~\ref{prop:basic} then says that the set 
$$E=\left\{p\in S\backslash (\mathrm{crit}(H_1)\cup \mathrm{crit}(H_2))\;|\;H_1(p)=H_2(p)\right\}$$
has an accumulation point in $S\backslash (\mathrm{crit}(H_1)\cup \mathrm{crit}(H_2))$.
According to Theorem~\ref{thm:Un}, we have $H_1=H_2$ 
on $S\backslash (\mathrm{crit}(H_1)\cup \mathrm{crit}(H_2))$.
By continuity it holds everywhere on $S$, and the conclusion follows.
\end{pf}

\section{Rational maps sharing a measure of maximal entropy}\label{sec:equilibrium}
\noindent It is an interesting problem to study rational maps that share a non-trivial Julia set, or have the same measure of maximal entropy.
Here we briefly summarize relevant notions and basic properties.
The reader can refer to \cite{DS} for details.

If $f$ is a rational map of $\hC$ of degree $d$, then the \textit{measure of maximal entropy} $\mu_f$ of $f$ is the unique
pull-back invariant probability measure with support equal to the Julia set $\J(f)$ of $f$,  achieving maximal entropy $\mathrm{log}\, d$. \emph{Pull-back invariant} 
means that $\frac{1}{d}f^*\mu_f=\mu_f$, or $\mu_f(f(A))=d\cdot  \mu_f(A)$, whenever $f$ is injective on a Borel set $A$. Note that this implies a weaker notion of invariance that
$f_*\mu_f=\mu_f$, or $\mu_f(f^{-1}(A))=\mu_f(A)$ for all Borel sets $A\sub \hC$.

For any probability measure $\mu$ on $\hC$ and a rational map $f$, 
the \textit{Jacobian} $J_\mu(f)$ of $f$ is defined as the Radon-Nikodym derivative of $f^*(\mu)$ with respect to $\mu$. Namely,
for each Borel measurable set $A$ such that $f$ is injective on $A$, 
$$\mu(f(A))=\int_A J_\mu(f)(x)\:\mathrm{d}\mu(x),\quad \mathrm{whenever}\; f|_A \;\mathrm{is}\; \mathrm{injective}.$$ 
If $\mu=\mu_f$ is the 
 measure of maximal entropy of a rational map $f$ of degree $\mathrm{deg}(f)$ at least 2, then $J_{\mu_f}(f)$ is constant equal to $\mathrm{deg}(f)$.
In fact, the converse is generally true (\cite[Theorem 1.118]{DS}):
\begin{thm}[\cite{DS}]\label{thm:constJac}
	Let $f$ be a ratinal map of degree $d$. If $\mu$ is a pull-back invariant probability measure of constant Jacobian $d$,
	then $\mu=\mu_f$.
\end{thm}
\begin{pf} The measure 
	$\mu_f$ is the unique probability measure such that the entropy $h_{\mu}(f)$ attains the maximum value $\mathrm{log}\, d$. The result now follows from the inequality~\cite{Pa} 
	$$h_\nu(f)\geq \int \mathrm{log}J_\nu(f)\;\mathrm{d}\nu,$$
	satisfied for any probability measure $\nu$ invariant under $f$, i.e., such that $f_*\nu=\nu$. Indeed, by setting $\nu=\mu$ from the statement of the theorem we get $h_\mu(f)\geq \mathrm{log}\, d$, and so $h_{\mu}(f)$ attains the maximum. Therefore $\mu=\mu_f$.  
\end{pf}  

It is well known that $\mu_f=\mu_{f^n}$ for all $n\in\N$, and commuting rational maps have the same measure of maximal entropy.
The following theorem is contained in~\cite{Le, Le2, LP}; see also~\cite[Theorem~1.7]{Ye}. 
\begin{thm}[\cite{Le, Le2, LP, Ye}]\label{thm:jutoeq}
Let $f, g\co\hC \to \hC$ be non-exceptional rational maps of degrees at least 2. 
Then $\mu_f=\mu_g$ if and only if for some $k, l\in \N$, 
 $$
  f^k\circ g^l=f^{2k}.
  $$
\end{thm}
  Here, a rational map is called \emph{exceptional} if it is conformally conjugate to either a power map, or a Chebyshev polynomial, or a Latt\`es map. In particular, a rational map whose Julia set is a Sierpi\'nski carpet is non-exceptional.
%
%
%
%

For general rational maps $f$ and $g$ sharing the same Julia set, they need not have common iterates up to  M\"obius transformations (\cite[Theorem 1.1]{Ye}). 
In general, holomorphic maps of the form $f^{-k}\circ f^k$ need not be M\"obius, as the following example shows. 
\begin{ex}\label{ex:simple}
  Consider $f(z)=z^2-\frac{1}{16z^2}$. It is postcritically finite and has a Sierpi\'nski carpet Julia set. A priori multi-valued map $f^{-1}\circ f(z)$ must be either $\pm z$ or $\pm \frac{i}{4z}$, i.e., a M\"obius transformation. 
  However, when $k\geq 2$, $f^{-k}\circ f^k$ first defined locally and then extended analytically may not be a M\"obius transformation. For example, by solving for the coefficients of a M\"obius transformation $\xi$ such that 
  \begin{equation}\label{eq:f2}
  	f^2\circ \xi=f^2,
  \end{equation}
  one can show that $\xi$ belongs to the group generated by $\rho(z)=\mathrm{i}z$ and $\iota(z)=\frac{1}{4z}$. 
  Indeed, if $\xi(z)={(az+b)}/{(cz+d)}$ satisfies~\eqref{eq:f2}, we may assume $a\neq 0$, otherwise we may consider $\iota\circ \xi$. Then 
  $$\lim_{z \to \infty}f^2(\xi(z))=\lim_{z \to \infty}f^2(z)=\infty,$$ while $\lim_{z \to \infty}\xi(z)=\frac{a}{c}$. Hence $f^2(\frac{a}{c})=\infty$, implying $c=0$. We may assume $d=1$, and thus $\xi(z)=az+b$. But $f^2(z)=z^4+O(1)$ and so $f^2(az+b)\thicksim (az+b)^4+O(1)$ as $z\to \infty$. These two functions being identical, it implies that $a^4=1$ and $b=0$ and the claim follows. In particular, the group of M\"obius transformation satisfuing~\eqref{eq:f2} has eight elements. On the other hand, the map $f^{-2}$ has 16 branches in a neighborhood of a noncritical value for $f^2$, so there must be some local 
  holomorphic maps of the form $f^{-2}\circ f^2$ which are not M\"obius.   
\end{ex}


\section{Proofs of Theorems~\ref{thm:main1} and~\ref{thm:main2}}\label{sec:main}
\noindent
We need the following lemma in the proof of Theorem~\ref{thm:main1}; cf.~\cite[Lemma~7.1]{BLM}.
\begin{lem} \label{lem:Rot} Let $\phi\co\partial \D \ra \partial \D$ be an orientation preserving branched covering, and suppose that there   exist numbers $k,l,n\in \N$, $k\ge 2$, such that 
\begin{equation} \label{eq:basiceq}
(P_l\circ \phi)(z)= (P_n\circ \phi\circ P_k)(z) \quad \text{for $z\in \partial \D$},
\end{equation}
where $P_n(z)=z^n$ for all $n\in\N$.
Then $l=n\cdot k$  and there exists  $b\in \partial \D, m\in \N$,  such that $\phi(z)=bz^m$ for all $z\in \partial \D$.
\end{lem}
\begin{pf} By counting topological degrees of both sides, one immediately sees that 
$l=n\cdot k$. 
Let $m$ be the degree of $\phi$. Then there exists a continuous function $\alpha\co \R \ra \R$ with $\alpha(t+2\pi)=\alpha(t)+2\pi m$ such that
$$ \phi(\mathrm{e}^{ \mathrm{i} t}) = \mathrm{e}^{ \mathrm{i}\alpha(t)} \quad \text{for  $t\in \R$}. $$
By  \eqref{eq:basiceq} we have 
$$ \mathrm{e}^{ \mathrm{i} l \alpha(t)} = \mathrm{e}^{ \mathrm{i} n\alpha(kt)}$$
for $t\in \R$. This implies that there exists a constant $c\in \R$ such that 
$$ \alpha(t)=\frac{n}{l}\alpha(kt)+c=\frac 1k \alpha(kt)+c  \quad \text{for  $t\in \R$}. $$
Using induction on $s\in\N\cup\{0\}$, one obtains
$$ \alpha\left(t+\frac{2\pi}{k^s}\right)=\frac 1k \alpha\left(kt+\frac{2\pi}{k^{s-1}}\right)+c=\frac 1k \alpha(kt)+\frac{2\pi}{k^s}m+c=\alpha(t)+\frac{2\pi}{k^s}m,$$  for all $t\in\R$.
Continuity of $\alpha$ then implies that $\alpha$ is a linear function and the conclusion follows.
\end{pf}

\subsection{Proof of Theorem~\ref{thm:main1}}
  We follow the steps, namely steps I-VI, in~\cite[Section 8]{BLM}.

 Steps I-IV are concerned with establishing equation
 \begin{equation}\label{eqn:geneq}
 g^{m'}\circ \xi = g^m\circ \xi \circ f^n,
 \end{equation}
 for some $n, m$, and $m'$, on all of $\J(f)$. 
 Proposition~\ref{prop:basic} and Theorem~\ref{thm:basicthm} above do just that under the assumption that $\xi\co \hC\to \hC$ is a quasiregular map with $\xi^{-1}(\J(g))=\J(f)$. 
 
 Following step V of~\cite[Section 8]{BLM}, we derive Equation~\eqref{eqn:fxi1}. Indeed, by post-composing both sides of Equation~\eqref{eqn:geneq} by a suitable iterate of $g$ and denoting $f^n$ by $F$ and $g^{m'-m}$ by $G$, we obtain
 \begin{equation}\label{eqn:speceq}
 G^{l+1}\circ \xi=G^l\circ\xi\circ F
 \end{equation}
 on $\J(f)$ for some $l\in\N$. This equation immediately implies Equation~\eqref{eqn:fxi1} on $\J(f)$.
 
 The final step VI of~\cite[Section 8]{BLM} defines a natural extension
 of $\xi$ from $\J(f)$ to $\hC$ that is analytic in each Fatou
 component of $\J(f)$. It is obtained by first extending $\xi|_{\J(f)}$
 in each periodic Fatou component of $\J(f)$ and then lifting via the
 dynamics to the remaining Fatou components, which are necessarily
 pre-periodic by Sullivan's No Wandering Domains Theorem. In~\cite{BLM}
 such extension is actually conformal in each Fatou component by
 \cite[Lemma~7.1]{BLM}. In our setting, i.e., when $\xi$ is assumed to be quasiregular, we
 obtain only an analytic extension to the periodic Fatou components,
 in a similar way using B\"ottcher coordinates and Lemma~\ref{lem:Rot}
 in place of~\cite[Lemma~7.1]{BLM}.
 
 Indeed, first let $U$ be a periodic Fatou component of $F$ from~\eqref{eqn:speceq}. We denote by $k\in\N$  the period of $U$, and define $V$ to be the Fatou component of $G$ bounded by $\xi(\partial U)$,  and $W=G^l(V)$, where $G$ and $l\in \N$ are as in~\eqref{eqn:speceq}.
 Then
 \eqref{eqn:speceq}  implies that $W$ is invariant under $G^k$.  By
 \cite[Lemma 3.3]{BLM}, the basepoint-preserving homeomorphisms  $\psi_U\colon (\overline U, p_U)\ra
 ( \overline \D, 0)$ and $\psi_W\colon(\overline W, p_W)\ra
 ( \overline \D,0)$, for unique $p_U\in U$ and $p_W\in W$,  conjugate $F^k$ and $G^k$, respectively, to power maps.
 Since $f$ and $g$ are postcritically finite, the periodic Fatou components $U$ and $W$ are superattracting, and thus the degrees of these power maps are at least 2.

 Again by Lemma~ \cite[Lemma 3.3]{BLM}, there exists $\psi_V\colon (\overline V, p_V)\to(\overline\D,0),\ p_V\in V$, such that
 $\psi_W\circ G^l\circ\psi^{-1}_V$ is a power map. Since $U,V,W$ are Jordan domains, the maps $\psi_U, \psi_V, \psi_W$ give  homeomorphisms between the boundaries of the corresponding Fatou components and $\partial \D$. Since $\xi$ is a
 branched covering of $\partial U$ onto $\partial V$, the map
 $\phi=\psi_V\circ \xi\circ\psi_U^{-1}$
 gives a branched covering  on $\partial \D$. Now \eqref{eqn:speceq} implies that on $\partial \D$ we have
 \begin{align*}
 	P_{d_3}\circ \phi&=\psi_W\circ G^{k+l}\circ \psi_V^{-1} \circ\phi
 	= \psi_W\circ G^{k+l}\circ \xi\circ \psi_U^{-1}\\
 	&=  \psi_W\circ G^{l}\circ \xi\circ F^k\circ \psi^{-1}_U
 	= \psi_W\circ G^{l}\circ \xi \circ \psi^{-1}_U \circ P_{d_1}\\
 	&=  \psi_W\circ G^{l}\circ  \psi^{-1}_V\circ \phi\circ P_{d_1}
 	= P_{d_2}\circ \phi \circ P_{d_1},
 \end{align*}
 for some $d_1, d _2, d_3\in\N$ with $d_1\geq 2$.  Lemma~\ref{lem:Rot} implies that $\phi$ extends to $\overline \D$ as a monomial around $0$, also denoted by $\phi$. In particular, $\phi(0)=0$, and so $\phi$ preserves the basepoint $0$ in $\overline \D$. If we  define $\xi=\psi_V^{-1}\circ \phi\circ\psi_U$  on $\overline U$, then $\xi$ is a branched covering 
 of $(\overline U,p_U)$ onto  $(\overline V, p_V)$.

 In this way,  we can analytically extend $\xi$ to every periodic Fatou component of $F$ so that $\xi$ maps the  basepoint of a Fatou component to the basepoint of the image component.
 To get an extension also for other, non-periodic, Fatou components of $F$, we proceed inductively
 on the level of such a given Fatou component $V$, i.e., the smallest $n\in\N\cup\{0\}$, such that $F^n(V)$ is periodic. Suppose that the level of $V$ is $n\in\N$ and that
 we already have an extension
 for all Fatou components with level smaller than $n$. This applies to the Fatou component
 $U=F(V)$ of $F$, and so a  conformal extension  $(U, p_U)\ra ( U',  p_{U'})$ of $\xi|_{\partial U}$ exists, where $U'$ is the Fatou component of $G$ bounded by $\xi(\partial U)$.
 Let $V'$ be the Fatou component of $G$ bounded by $\xi(\partial V)$, and $W=G^{l+1}(V')$.
 Then, by using \eqref{eqn:speceq} on $\partial V$ we conclude that $G^l(U')=W$.
 Define $\alpha=G^l\circ \xi|_{\overline U}\circ F|_{\overline V}$ and
 $\beta=\xi|_{\partial V}$. Then the assumptions of the Lifting Lemma~\cite[Lemma 3.4]{BLM} are satisfied for $D=V$, $p_D=p_{V}$, and the iterate
 $G^{l+1}\co V'\ra W$ of $G$. Indeed, $\alpha$ is continuous on $\overline V$ and holomorphic
 on $V$, we have
 $$\alpha^{-1}(p_W)=F|_{\overline V}^{-1}(  \xi|_{\overline U}^{-1}( p_{U'}))=F|_{\overline V}^{-1}(p_U)=p_V,$$
 and
 $$G^{l+1}\circ \beta= G^{l+1}\circ  \xi|_{\partial V}= G^{l}\circ \xi|_{\partial U}\circ F|_{\partial V}
 =\alpha|_{\partial V}.$$ Since $\beta$ is continuous, \cite[Lemma 3.4]{BLM} implies that there exists a continuous map $\widetilde \alpha$  of $(\overline V, p_V)$ onto $(\overline V', p_{V'})$ such that $\widetilde \alpha|_{\partial V}=\beta=\xi|_{\partial V}$. In other words,
 $\widetilde \alpha$ gives the desired basepoint preserving analytic extension to the Fatou component $V$.

  
  Now, $\xi$ has a necessarily unique extension $\widehat{\xi}$, analytic in each Fatou component of $F$ and hence of $f$. It is not hard to see, arguing as in~\cite[Lemma~2.1]{BLM}, 
  that since $\xi$ is quasiregular, so is $\widehat{\xi}$. Indeed, the arguments in the proof of~\cite[Lemma~2.1]{BLM} can be localized by first observing that there are only finitely many Fatou components $U$ of $f$ such that $\xi|_{\partial U}$ is not a homeomorphism, and then using the fact that the boundaries of Fatou components of $f$ are quasicircles, and hence removable for quasiregular maps.  Sto\"ilow's Factorization Theorem now tells us that 
  $\widehat{\xi}=h\circ q$ for some quasiconformal map $q$ on $\hC$ and a holomorphic map $h$ on $q(\hC)$.
  It follows that $q$ is conformal on the full measure set $\hC\backslash \J(f)$, and hence a M\"obius transformation. 
  Then $h$ is a holomorphic map on $q(\hC)=\hC$, i.e., $h$ is rational, and so is $\widehat{\xi}$. 
  
  Finally, if we replace $\xi$ by $\widehat\xi$, both sides of Equation~\eqref{eqn:speceq} are rational maps that agree on $\J(f)$. Therefore they agree on all of $\hC$, establishing Equation~\eqref{eqn:fxi1}.
\qed


\subsection{Proof of Theorem~\ref{thm:main2}}
The measure of maximal entropy  $\mu_{\widehat\xi}$ for $\widehat\xi$ is a probability measure supported on $\J(\widehat\xi)$. By passing to iterates and noting that the Julia sets do not change under taking iterates, we may assume that Equation~\eqref{eqn:fxi1} is satisfied for $f$ and $g$ rather than $F$ and $G$, respectively. Choosing $f=g$, it gives that $f^n\circ \widehat\xi$ commutes with $f^k$ for some $k,n\in\N$. 
Since we can choose arbitrary $k$ in Equation~\eqref{eqn:fxi1} and also post-compose both sides of this equation by an arbitrary iterate of $f$, we can assume that $n=k$, i.e., $f^k\circ \widehat\xi$ commutes with $f^k$ for some $k\in\N$. 
Therefore, $f^k\circ \widehat\xi$ and $f^k$
share a measure of maximal entropy  $\mu=\mu_f=\mu_{f^k}=\mu_{f^k\circ \widehat{\xi}}$. Calculating the Jacobian of $\widehat{\xi}$ with 
respect to $\mu$, we obtain 
$$J_\mu(\widehat{\xi})(x)=J_\mu(f^k\circ \widehat{\xi}(x))/J_\mu(f^k)(\widehat{\xi}(x))=\mathrm{const}/\mathrm{const}=\mathrm{const},$$ i.e., $(\widehat{\xi})^*(\mu)=\mathrm{deg}(\widehat\xi)\mu$. 
By Theorem~\ref{thm:constJac}, this implies $\mu=\mu_{\widehat{\xi}}$, and in particular $\J(\widehat{\xi})=\J(f^k)=\J(f)$. 
Theorem~\ref{thm:jutoeq} applied to $f$ and $\widehat\xi$ gives Equation~\eqref{eqn:fxi2}.  Equation~\eqref{eqn:fxi2} in turn gives that the critical points of $\widehat\xi$ are among the critical points of an iterate of $f$, and indeed $\widehat\xi$ is postcritically finite.
\qed

\bigskip

Several comments are in order. If $\xi$ is a quasiconformal map of $\widehat\C$ in Theorem~\ref{thm:main2}, then~\cite[Theorem~1.1]{BLM} gives that $\xi$ is the restriction of a M\"obius transformation $\widehat\xi$ to $\J(f)$.
From \cite[Theorem 1.1]{Le} we also know that  then $\widehat{\xi}^m=\mathrm{id}$ for some $m\in\N$, and the arguments from~\cite{BLM} give that
$\widehat{\xi}$ satisfies the functional equation 
$$
f^{2k}\circ \widehat{\xi} = f^k\circ \widehat{\xi} \circ f^k \;\mathrm{for\ some}\; k\in \N;
$$
see also the proof of Theorem~\ref{thm:main2} above. Replacing $f$ by its iterate, we can assume $f^2\circ \widehat{\xi}=f\circ \widehat{\xi}\circ f$. 
\begin{question}
  Let $f$ be a postcritically finite rational maps whose Julia set $\J(f)$ is a Siepi\'nski carpet, and let $\widehat\xi$ be a M\"obius transformation preserving $\J(f)$, and such that $f^2\circ \widehat{\xi}=f\circ \widehat{\xi}\circ f$. Does it imply that $F\circ \widehat{\xi}=F$, where $F$ is an iterate of $f$?
\end{question}
In fact, it is unclear whether such relations hold for specific maps such as  $f(z)=z^2-\frac{1}{16z^2}$. However, we suspect that the following question has a positive answer.
\begin{question}
Let $f(z)=z^2-\frac{1}{16z^2}$, and  $\rho,\iota$ be as in the example of Section~\ref{sec:equilibrium}, i.e., $\rho(z)=\mathrm{i}z$ and $\iota(z)=\frac{1}{4z}$. 
Is it true that the group of all M\"obius transformations preserving $\J(f)$ is $\left\langle \rho,\iota\right\rangle$?
\end{question}

\section{Quasiregular maps between other topological types of Julia sets}\label{sec:other}
\noindent 
One can ask whether similar statements to those in Theorems~\ref{thm:main1} and Theorem~\ref{thm:main2} hold for other topological types of non-exceptional postcritically finite rational maps. Here, again, by exceptional maps we mean, up to conjugation, power maps, Chebyshev polynomials, or Latt\`es maps. 
In this section we discuss two examples: the basilica and the {Apollonian gasket} cases. In contrast to the Siepi\'nski carpet case, we demonstrate that a quasiregular map of the basilica or the gasket Julia sets need not be related to the dynamics of the corresponding map, contrasting the Sierpi\'nski carpet case of Theorems~\ref{thm:main1} and~\ref{thm:main2}.

\begin{ex} The basilica; see~\cite{LM} for a description of the quasisymmetry group of the basilica and the discussion of the Thompson group $T$, and~\cite{BF} for more general statements and results.
The \textit{basilica} is the filled Julia set $\mathcal{K}(f)$ of the quadratic polynomial $f(z)=z^2-1$.
This polynomial has a superattracting cycle $\gamma=\{0,1\}$. By definition, the \emph{immediate basin} of $\gamma$ is the union of two components $U_0$
and $U_{-1}$ of the Fatou set containing 0 and -1, respectively. These components are Jordan disks, and they have a unique intersection point 
$\alpha=\frac{1-\sqrt{5}}{2}$, which is a fixed point of $f$. 

Let $U_\infty$ denote the basin at infinity, i.e., the unbounded Fatou component of $f$, and let $\phi_\infty$ be the B\"ottcher coordinate
of $U_\infty$. Namely, $\phi_\infty$ is the conformal map of $U_\infty$ onto $\C\backslash \overline{\D}$ conjugating $f$ to $g(z)=z^2$.
It is known that the Julia set of a hyperbolic map $f$ is locally connected, so $\psi_\infty=\phi_\infty^{-1}\co \C\backslash \overline{\D}\to U_\infty$ exteneds continuously to 
the boundary, and so induces a continuous boundary map $\psi$ of $\mathbb{S}^1$ onto the Julia set $\J(f)$ of $f$ satisfying $\psi\circ g=f\circ \psi$.
Since $\alpha$ is a fixed point of $f$, it can be calculated that 
$\psi^{-1}(\alpha)=\{\mathrm{e}^{\frac{2}{3}\pi \mathrm{i}},\mathrm{e}^{\frac{4}{3}\pi \mathrm{i}}\}$. 

Since the interior of $\mathcal{K}(f)$ for $f(z)=z^2-1$ equals the (full, not just immediate) basin of the cycle $\gamma$, every bounded Fatou component $U$ of $f$ eventually, i.e., under a certain iterate of $f$, lands in the cycle $\{U_0, U_{-1}\}$. 
In fact, for each such bounded Fatou component $U\neq U_0$ there exists a unique $n\in\N$ such that $f^n\co U\to U_0$ is a conformal map. Let $L(U)$ denote the limb of $U$, namely the closure of the union of $U$ and all the components of the interior of $\mathcal{K}(f)$ separated by $U$ from $U_0$. Note here that the ba\-si\-li\-ca has a tree-like structure.  Then $\psi^{-1}(\partial L(U))$ maps under $g^n$ onto 
$[{e}^{-\frac{2}{3}\pi {i}}, {e}^{\frac{2}{3}\pi {i}}]\sub \mathbb{S}^1$. The two endpoints ${e}^{\pm\frac{2}{3}\pi {i}}$ are mapped to the same point under $\psi$, called the \emph{root} of $L(U)$. 

In~\cite[Lemma~5.4]{LM}, it is shown that the Thompson group $T$ acts by quasisymmetries on $\J(f)$, where $f(z)=z^2-1$, leaving $\partial U_0$ invariant. In fact, such an action is defined piecewise using forward iterates of $f^2$ and branches of $f^{-2n},\ n\in\N$. Thompson group elements shuffle the limbs rooted at $\partial U_0$. Moreover, these quasisymmetries extend to global quasiconformal maps of $\C$. It is immediate that if one takes a non-trivial element of the Thompson group acting by a quasisymmetry $\xi$ on $\J(f)$, say that fixes all the limbs rooted on a certain non-trivial arc of $\partial U_0$, then this $\xi$ does not relate globally to the dynamics of $f$, say as in Equation~\eqref{eqn:fxi2}. Such a map $\xi$ has degree one however. Post-composing $\xi$ with $f$ produces a quasiregular map of degree two that does not satisfy Equation~\eqref{eqn:fxi2} on $\J(f)$.   


\end{ex}

\begin{ex} The Apollonian gasket; see~\cite{LLMM} for the discussion of the quasisymmetry group of the Apollonian gasket Julia set, and also \cite{LMMN} for more general statements.

 We need the following quasiregular extension lemma; see, e.g., \cite[Lemma~2.2]{LN} for similar results.
\begin{lem}\label{lem:qrext}
	Let $f$ be a continuous orientation preserving  branched covering map of the unit circle $\mathbb{S}^1$, and for $n\in\N$, let $I_1, I_2,\dots, I_n$, be closed arcs of $\mathbb{S}^1$ with disjoint interiors and whose union is $\mathbb{S}^1$.
	Assume that $f|_{I_j}$ extends in a neighborhood of $I_j$ to a holomorphic function with non-vanishing derivative for each $j=1,2,\dots,n$.
%
 Then  $f$ can be extended to a quasiregular map $\widehat{f}\co\overline{\D}\to \overline{\D}$.
\end{lem}
\begin{pf}
Let $e^{i\theta_j},\ j=1,2,\dots,n$, be the  points on $\mathbb{S}^1$ such that $e^{i\theta_{j}}\in I_{j}\cap I_{j+1},\ j=1,2,\dots n-1$, and $e^{i\theta_n}\in I_n\cap I_1$. We define $\widehat{f}(z)=|z|f(\frac{z}{|z|})$ on $\overline{\D}\setminus\{0\}$, and $\widehat{f}(0)=0$. It is a continuous map extending $f$. If $z\in \D$ and ${\rm arg} (z)\notin\{\theta_1, \theta_2,\dots, \theta_n\}$, where ${\rm arg} (z)$ denotes the argument of $z$, we have
  $$\partial_{\bar{z}}\widehat{f}=\frac{z}{2|z|}f\left(\frac{z}{|z|}\right)-\frac{1}{2}\left(\frac{z}{|z|}\right)^2 f'\left(\frac{z}{|z|}\right),$$
  $$\partial_z\widehat{f}=\frac{\bar{z}}{2|z|}f\left(\frac{z}{|z|}\right)+\frac{1}{2}f'\left(\frac{z}{|z|}\right).$$
  Since $f$ maps $\mathbb{S}^1$ to itself and preserves the orientation, 
we get that $\mathrm{arg}(f'(z))+\mathrm{arg}(z)\equiv \mathrm{arg}(f(z))\; \mathrm{mod}\;2\pi$ for each $z\in \mathbb{S}^1,\ {\rm arg} (z)\notin\{\theta_1, \theta_2,\dots, \theta_n\}$. 
  Therefore, for $z\in\D,\ {\rm arg} (z)\notin\{\theta_1, \theta_2,\dots, \theta_n\}$,  $$\left\lvert \partial_{\bar{z}}\widehat{f}\right\rvert =\frac{1}{2}\left\lvert 1-\left\lvert f'\left(\frac{z}{|z|}\right) \right\rvert\right\rvert,\quad
  \left\lvert \partial_z\widehat{f}\right\rvert =\frac{1}{2}\left(1+\left\lvert f'\left(\frac{z}{|z|}\right) \right\rvert\right).$$
  By the assumption, there exists $\varepsilon>0$ such that $|f'(z)|\geq \varepsilon$ for any $z\in \mathbb{S}^1,\ {\rm arg} (z)\notin\{\theta_1, \theta_2,\dots, \theta_n\}$, and hence  
$$\left\lvert\frac{\partial_{\bar{z}}\widehat{f}}{\partial_z\widehat{f}}\right\rvert=\frac{\left\lvert 1-\left\lvert f'\left(\frac{z}{|z|}\right) \right\rvert\right\rvert}{1+\left\lvert f'\left(\frac{z}{|z|}\right) \right\rvert}\leq k<1,$$ 
for some $k$ independent of $z$. From the assumptions and the extension formula it is immediate that $\widehat f$ is absolutely continuous on lines, and hence $\widehat{f}$ is a quasiregular extension of $f$.  
\end{pf}

Now we are able to prove the following result.
\begin{thm}
  There exists a postcritically finite rational map $f\co \mathbb{C} \to \mathbb{C}$ with gasket Julia set $\J(f)$, 
  and a global quasiregular map $\xi\co \hC \ra \hC$ with $\xi^{-1}(\J(f))=\J(f)$, 
  such that $\xi|_{\J(f)}$ does not have a rational extension to $\hC$. 
\end{thm}

\begin{pf}
  Let $f(z)=\frac{3z^2}{2z^3+1}$ considered in~\cite{LLMM}. 
  Its critical points are $0,1,\omega$ and $\omega^2$, where $\omega=e^{2\pi i/3}$, with orbits $0\mapsto 0$, $1\mapsto 1$, and $\omega\mapsto \omega^2\mapsto \omega$. In particular, $f$ is postcritically finite.  Note that $\bar{f}$, where the bar denotes the complex conjugation, is critically fixed, i.e., all its critical points are fixed. The Julia set $\J(f)=\J(\bar{f})$ is the Apollonian gasket depicted in Figure~\ref{fig:gasket2}. 
  
 \begin{figure}[h]
 	\centering
 	\includegraphics[width=0.4\textwidth]{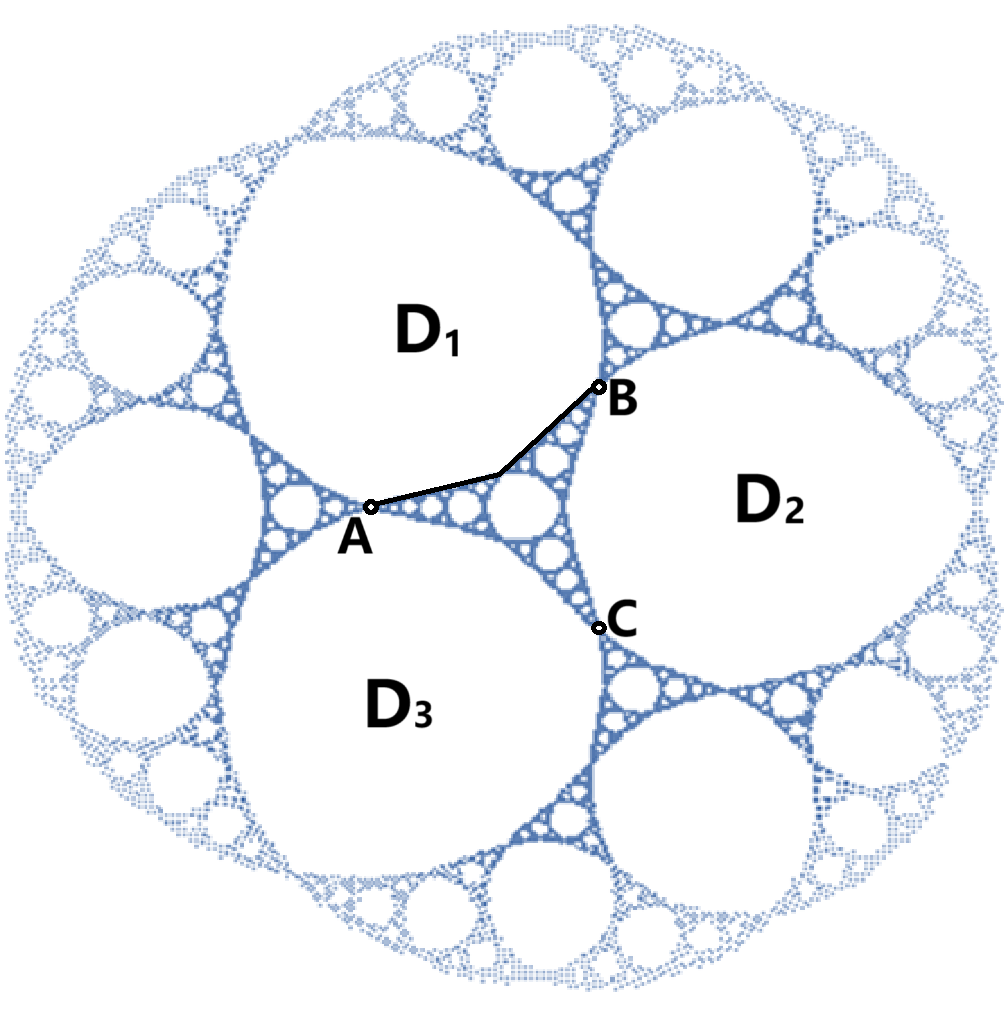}
 	\caption{The Julia set of
 		$f(z)=\frac{3z^2}{2z^3+1}$. 
 		}
 		\label{fig:gasket2}
 \end{figure}

  Denote the central closed triangular region in Figure~\ref{fig:gasket2} by $\Delta ABC$, with three edge $AB, BC, CA$ being arcs on the boundaries of the disks $D_1, D_2, D_3$, respectively. The map $f^2$ fixes the points $A, B, C$, and $D_1,D_2,D_3$,  are its fixed Fatou components.
  Now we take 
  $$\xi(z)= \begin{cases}
  f^2(z),  \quad z\in \Delta ABC  \\
  z,\quad z\in \hC\backslash \left((\Delta ABC)\cup \sqcup_{i=1}^{3}D_i\right).
  \end{cases}
  $$
  The map $\xi$ can be extended to a quasiregular map of $\hC$. Indeed, the boundaries of $D_1, D_2, D_3$ being quasicircles, it is enough to show that $\xi$ can be extended as a quasiregular map in each of these domains.
  We extend it to $D_1$ and the other two are done in the same way. 
  Let $\phi\co\overline{D_1}\to \overline{\D}$ be the B\"ottcher coordinate of $\overline{D_1}$.
  Counting topological degrees, for example, we find that 
  $$\phi\circ f^2=P_4\circ \phi \quad \mathrm{on} \quad \overline{D_1},$$ where $P_d(z)=z^d$ as above. Therefore, $\xi\co AB\to \partial D_1$ is conjugated to 
  $$P_4(z)=z^4\co\phi(AB)\to \mathbb{S}^1.$$
  Applying Lemma~\ref{lem:qrext} to the conjugate of $\xi$ by $\phi$, we obtain a quasi\-re\-gu\-lar extension into the unit disk. Conjugating it back by $\phi$, which is a conformal map, we obtain a quasiregular extension of $\xi$ into $D_1$, and hence a quasiregular extension of $\xi$ to $\hC$, still denoted by $\xi$. 
  
  Finally, $\xi|_{\J(f)}$ cannot be extended to a global rational map. Indeed, if it extended to a rational map of $\widehat{\C}$, then, by uniqueness for analytic maps, it would have to equal to $f^2$ in $\Delta ABC$ and to the identity in $\hC\backslash \left((\Delta ABC)\cup \sqcup_{i=1}^{3}D_i\right)$. This is a contradiction since $f^2$ is not equal to the identity.    
\end{pf}

\end{ex}



\begin{thebibliography}{BHa}


\bibitem[AIM]{AIM} K.~Astala, T.~Iwaniec, G.~M.~Martin, \emph{Elliptic partial differential equations and quasiconformal mappings in the plane}, Princeton Univ.\ Press, Princeton, NJ, 2009.  


\bibitem[Be]{Be} A. F. Beardon, {\em Iteration of rational functions}, Springer, New York, 1991. 

\bibitem[BF]{BF} J.~Belk, B.~ Forrest,  \emph{Quasisymmetries of finitely ramified Julia sets}, Math.\ Ann.\ 393 (2025), no.\ 2, 1683--1740.


\bibitem[Bo]{Bo} M. Bonk, \emph{Uniformization of Sierpi\'nski carpets in the plane} Invent. Math., 186 (2011), 559--665.


\bibitem[BKM]{BKM} M.~Bonk, B.~Kleiner, S.~Merenkov, \emph{Rigidity of Schottky sets},
Amer.\ J.\ Math.\ 131 (2009), no.\ 2, 409--443.


\bibitem[BLM]{BLM} M. Bonk, M. Lyubich, S. Merenkov, \emph{Quasisymmetries of Sierpi\'nski carpet Julia sets} Adv. Math, 301 (2016), 383--422.



\bibitem[DS]{DS} T. Dinh, N. Sibony, {\em Dynamics in Several Complex Variables: Endomorphisms of Projectice Spaces and Polynomial-like Mappings}, Lecture Notes in Mathematics, Springer, 1998.


\bibitem[Le]{Le} G. M. Levin, \emph{Symmetries on Julia sets} (Russian), Mat. Zametki 48 (1990),   72--79, 159; translation in Math. Notes 48 (1990),  1126--1131 (1991).

\bibitem[Le2]{Le2} G. M. Levin, \emph{Letter to the Editor}, Math. Notes 69 (3) (2001), 432--433.


\bibitem[LP]{LP} G.~Levin, F.~Przytycki, 
\emph{When do two rational functions have the same Julia set?}
Proc. Amer. Math. Soc. 125 (1997), no. 7, 2179--2190.


\bibitem[LLMM]{LLMM} R.~Lodge, M.~Lyubich, S.~Merenkov, S.~Mukherjee,  \emph{On dynamical gaskets generated by rational maps, Kleinian groups, and Schwarz reflections} 
Conform. Geom. Dyn. 27 (2023), 1--54.

\bibitem[LP2]{LP2} R. Luisto, P.  Pankka, \emph{Sto\"ilow’s theorem revisited} Expo. Math, 38 (2020), 303--318.

\bibitem[LN]{LN} Y. Luo, D. Ntalampekos, \emph{Uniformization of gasket Julia sets} Preprint, arXiv:2411.17227.

\bibitem[LM]{LM} M. Lyubich, S. Merenkov, \emph{Quasisymmetries of the basilica and the Thompson group}, Geom. Funct. Anal. 28 (2018) 727--754.

\bibitem[LMMN]{LMMN} M.~Lyubich, S.~Merenkov, S.~Mukherjee, D.~Ntalampekos, \emph{David extension of circle homeomorphisms, welding, mating, and removability}, Mem.\ Amer.\ Math.\ Soc.\ 313 (2025), no.\ 1588, v+110 pp.
   


\bibitem[Me1]{Me1} S.~Merenkov, \emph{Planar relative Schottky sets and quasisymmetric maps},
Proc.\ Lond.\ Math.\ Soc.\ (3) 104 (2012), 455--485.


\bibitem[Me2]{Me2} S.~Merenkov,  \emph{Local rigidity of Schottky maps}, Proc.\ Amer.\ Math.\ Soc.\ 142 (2014), no. 12, 4321--4332.


\bibitem[Me3]{Me3} S.~Merenkov,  \emph{Local rigidity for hyperbolic groups with Sierpi\'nski carpet boundaries} Compo. Math. 150 (2014), 1928--1938.


 
\bibitem[Pa]{Pa} W.~Parry, \emph{Entropy and generators in ergodic theory}, W.\ A. Benjamin, Inc., New
York-Amsterdam, 1969. 
 
\bibitem[St]{St} S.~Sto\"ilow, \emph{Sur les transformations continues et la topologie des fonctions analytiques}, Ann. Sci.\ \'Ec.\ Norm.\ Sup\'er. (3) 45 (1928) 347--382.


\bibitem[Wh]{Wh} G. T. Whyburn, \emph{Topological characterization of the Sierpi\'nski curve}, Fund. Math. 45 (1958), 320--324.




\bibitem[Ye]{Ye} H.~Ye, \emph{Rational functions with identical measure of maximal entropy},
Adv. Math. 268 (2015), 373--395.


\end{thebibliography}
\end{document}